\documentclass[11pt]{amsart}
\usepackage{amsmath
,amscd,graphicx
}
\makeatletter
\atdef@ N#1N#2N{\CD@check{N..N..N}{\llap{$\m@th\vcenter{\hbox
  {$\scriptstyle#1$}}$}\phantom{\Big\downarrow} % \Big\downarrow
  \rlap{$\m@th\vcenter{\hbox{$\scriptstyle#2$}}$}&&}}
\makeatother
\theoremstyle{plain}
\newtheorem{theorem}{Theorem}[section]
\newtheorem{lemma}[theorem]{Lemma}
\newtheorem{corollary}[theorem]{Corollary}

\newtheorem{ther}{Theorem}

\newtheorem*{claim}{Claim}
\newtheorem*{acknowledgement}{Acknowledgements}

\theoremstyle{remark}
\newtheorem*{remark}{Remark}
\newtheorem*{notation}{Notation}
\newtheorem*{example}{Example}

\theoremstyle{definition}
\newtheorem{defn}[theorem]{Definition}
\newtheorem{convention}[theorem]{Convention}
\newtheorem{construction}[theorem]{Construction}

\def \R {\mathbf{R}}
\def \Z {\mathbf{Z}}
\def \C {\mathbf{C}}
\def\RP{\mathbf{RP}}
\def \HH{\mathcal{H}}
\def\zz{\ifmmode{z}\else{$z$}\fi}
\def\abar{\bar{a}}
\def\bbar{\bar{b}}
\def\mtau{M^\tau}
\def\wtau{W^\tau}
\def\ytau{Y^\tau}
\def\M{\mathcal{M}}
\def\hflat{$h$-flat\ }
\def\MZ{\mathcal{M}^{\mathbf{Z}_2}}
\def\MU{\mathcal{M}^{\mathrm{U}(1)}}
\def\MSU{\mathcal{M}^{\mathrm{SU}(2)}}
\def\mapt{S^1 \times_\tau \Sigma}
\def\mapphi{S^1 \times_\varphi \Sigma}
\DeclareMathOperator{\SU}{SU}
\DeclareMathOperator{\SF}{SF}
\DeclareMathOperator{\SO}{SO}
\DeclareMathOperator\ad{ad}
\def\ada{\ad\alpha}
\def\su{\ifmmode{\SU(2)}\else{$\SU(2)$}\fi}
\def\trans{\cap\hskip -1.62ex |\hskip 1.1ex} 
\DeclareMathOperator{\U}{U}
\DeclareMathOperator{\SL}{SL}
\DeclareMathOperator{\CS}{CS}

\DeclareMathOperator{\CF}{CF}
\DeclareMathOperator{\HF}{HF}

\DeclareMathOperator{\rep}{Rep}
\begin{document}

\baselineskip.525cm
\title{Mutation and Gauge theory I: Yang-Mills Invariants}
\thanks{The author was partially supported by NSF Grant 4-50645}
\author[Daniel Ruberman]{Daniel Ruberman}
\address{Department of Mathematics\newline
Brandeis University \newline
Waltham, MA 02254}
\email{\rm{ruberman@binah.cc.brandeis.edu}}
\date\today
\maketitle

Mutation is an operation on 3-manifolds containing an embedded surface $\Sigma$
 of 
genus $2$.  It is defined using the unique involution, called $\tau$, of
$\Sigma$ with the property that $\Sigma/\tau \cong S^2$.   In brief, given a
$3$-manifold $M$ containing $\Sigma$, its mutant $\mtau$ is obtained by
cutting $M$ along $\Sigma$, and regluing using $\tau$.  This operation was
introduced in~\cite{ruberman:mutation} as the analogue for closed manifolds of
the mutation operation on knots described by Conway~\cite{conway}.
It is not easy to distinguish a $3$-manifold from its mutant---there is a long list of
invariants which they have in common:  their Gromov norm~\cite{ruberman:mutation},
Reidemeister torsion~\cite{porti},
Chern-Simons and
$\eta$-invariants~\cite{meyerhoff-ruberman:I,meyerhoff-ruberman:II}
(if $M$ is hyperbolic),  Casson's invariant~\cite{kirk:casson}, and many of 
the `quantum' invariants of
3-manifolds~\cite{kania,kawauchi:quantum,lickorish:same_su2,RGK:mutation,
rong:mutation}.

In this article, we will show that the 
instanton Floer homology~\cite{floer:instanton} and 
$\Z$-graded instanton homology~\cite{fs:graded} of homology spheres are unchanged by
mutation.  
\begin{ther}\label{hfmut} Let $M$ be an oriented homology $3$-sphere, with 
(instanton) Floer homology $HF_*(M)$, which contains a genus-$2$ surface,
and let $\mtau$ be the result of mutation along $\Sigma$.
Then $HF_*(M) \cong HF_*(\mtau).$  Similarly, if $HF^\mu$ denotes the $\Z$-graded
instanton homology of Fintushel-Stern, then $HF^\mu_*(M) \cong HF^\mu_*(\mtau).$ 
\end{ther}

In section~\ref{4mfd} we will define two types of mutation operations on certain
$4$-manifolds, and show that they preserve the Donaldson invariant.  A companion article
(in preparation) will show that the 
$3$-dimensional Seiberg-Witten analogue of Casson's invariant is unchanged by mutation. 

Theorem~\ref{hfmut} provides an alternate
proof of P.~Kirk's result~\cite{kirk:casson} on the Casson invariant.  In his
paper~\cite{kawauchi:mutative}, A.~Kawauchi constructs homology 
spheres with the same Floer homology as their mutants, and remarks that the
general case does not seem to be known.  On reviewing Kawauchi's paper 
for Mathematical
Reviews, my interest in the problem was stimulated by this remark.  The
papers of Kirk and Kawauchi are based on the connection between surgery and Casson's
invariant (respectively Floer homology), whereas we will proceed directly from 
the definition of Floer homology in terms of \su-representations.
The restriction to homology spheres is largely for technical convenience; it is likely 
that the
proof of theorem~\ref{hfmut} would presumably extend to Floer-type
theories~\cite{austin-braam:floer,braam-donaldson:gluing} defined for more general
$3$-manifolds.
\begin{acknowledgement}
I would like to thank Tom Mrowka for several helpful conversations, and Hans Boden for
pointing out the example discussed after theorem~\ref{donaldson}. 
\end{acknowledgement}

\begin{notation}
For the rest of the paper, $M$ will denote a closed, oriented $3$-manifold,
and $\Sigma$ a genus-$2$ surface.  For any space $X$, and Lie group $G$,
the space of representations of $\pi_1(X)$ into $G$ will be denoted $\rep(X,G)$.
The equivalence  classes (under the relation of conjugacy by elements of $G$) will
be denoted $\chi(X,G)$.  If it is obvious what group is being discussed, then
the `$G$' may be dropped.
\end{notation}

\section{Mutation and Floer Homology}
The involution $\tau$ has several related properties which are
responsible for the equality of the invariants cited above after mutation.
The basic one is that any simple closed curve $\gamma$ on $\Sigma$ is isotopic
to one which is taken
to itself by $\tau$, perhaps with a reversal of orientation.  This implies that
$\tau$ is in the center of the mapping class group, but more importantly for
our purposes, implies the following lemma, well-known in certain circles.
\begin{lemma}\label{symmrep}  Let $\varphi:\pi_1(\Sigma)\to \su$ be a
representation.  Then
$\varphi\circ\tau_*$ is conjugate to $\varphi$.  The same is true for representations
of $\pi_1(\Sigma)$ in $\SL_2(\R)$.
\end{lemma}
A similar lemma, concerning instead representations in $\SL_2(\C)$,
may be found in \S2 of~\cite{ruberman:mutation}, and the proof there may be
adapted, {\sl mutatis mutandis}.    
\begin{remark}  
In the field of `quantum invariants' of $3$-manifolds, this lemma is seen 
as a reflection of a self-duality of certain representations
of $\su$.  This duality does not hold for larger rank unitary
groups, so one expects that invariants based on, say, $\SU(3)$
representations, would change under mutations.  We will return to this
point in Section~\ref{4mfd}.
\end{remark}

Using the standard correspondence between representations and flat connections, this
lemma says that for any flat $\su$-connection $\alpha$ on $\Sigma$, the pull-back
$\tau^*\alpha$ is gauge equivalent to $\alpha$. 
If $\alpha $ is reducible, then there is a choice of (constant) gauge
transformation $g \in \rm{Stab}(\alpha)/\Z_2$ with $\tau^*\alpha= g^*\alpha$.  Any
such $g$ defines an automorphism $\hat\tau$ of the connection $\alpha$ on the trivial
bundle, as the following composition:
\begin{equation}\label{taut}
\begin{CD}
\Sigma \times \su  @>g^{-1}>> \Sigma \times \su  @>\tau \times \mathrm{id}>>
\Sigma \times \su \\
@VVV @VVV @VVV\\
\Sigma  @>\mathrm{id}>> \Sigma @>\tau>> \Sigma
\end{CD}
\end{equation}
By definition, $\hat\tau$ covers $\tau$, and induces an automorphism $\tau^*$ of the
$su_2$-valued forms $\Omega^*(\Sigma;\ada)$, and the twisted cohomology groups
$H^*(\Sigma;\ada)$.

Starting from Lemma~\ref{symmrep}, there is an obvious path to take in showing
the equality of $HF_*$ of mutant homology spheres.  Let $\CF_*$ denote the chain
complex which computes the instanton homology; in favorable circumstances this
has a basis  indexed by  the flat $\su$-connections on $M \times \su$, or
equivalently by $\chi(M,\su)$.  The $\Z/8\Z$-grading is given by spectral flow.
For $\mu\in \R$ such  that $\mu\neq \CS(\alpha)$ for any $\alpha \in \chi(M)$,
Fintushel and Stern~\cite{fs:graded} have defined $\Z$-graded chain groups
$\CF^\mu_*$.  These have the same basis as $\CF_*$, but the grading is lifted from 
$\Z/8\Z$ to $\Z$ using the monotonicity properties of the Chern-Simons invariant. 

Suppose that $\Sigma \subset M$, separating $M$ into
two components whose closures will be denoted $A$ and $B$.  When $M$ has a Riemannian
metric, the metric will be assumed to be a product in a neighborhood of $\Sigma$. 
Moreover, we will assume that $\tau$ is an isometry of the restriction of the metric to
$\Sigma$.   In this notation,
$\mtau$ will be given by $A\cup_\tau B$, and inherits a metric from $M$.  Given a 
representation $\varphi$ of $\pi_1(M)$ in $\su$,
let $\varphi_\Sigma$ denote its restriction to
$\pi_1(\Sigma)$, with $\varphi_A$ (resp.~$\varphi_B$) the restrictions to
$\pi_1(A)$ (resp.~$\pi_1(B)$.)  Choose an element $g \in \su$
with $\varphi_\Sigma \circ \tau_* = g^{-1}\varphi_\Sigma g$, and conjugate
$\varphi_B$ by $g$, to get a representation $\varphi^\tau$ of 
$\pi_1(\mtau)$.  Thus we get a sort of correspondence between
$\chi(M,\su)$ and $\chi(\mtau,\su)$, which should lead to an isomorphism on
instanton homology.  

There are several issues with which one must deal:
\begin{itemize}
\item[(i)] The character variety $\chi(M)$ may not consist of a finite number of 
smooth points.
\item[(ii)] If $\varphi_\Sigma$ is reducible, then there is a choice (parameterized
by $\rm{Stab}(\varphi_\Sigma)/\Z_2$) of elements $g$ conjugating $\varphi_\Sigma$ to
$\varphi_\Sigma \circ \tau_*$.
\end{itemize}
Even if these problems do not arise, so that we have a sensible map $T_*:
\CF_*(M) \to \CF_*(\mtau)$ (and $T^\mu_*$ in the $\Z$-graded case) we would need to show:
\begin{itemize}
\item[(iii)] The map $T_*$ preserves the $\Z/8\Z$ grading in $\CF_*$ and $T^\mu_*$ the $\Z$
grading in $\CF^\mu_*$ .
\item[(iv)] $T_*$ and $T^\mu_*$ are  chain maps.
\end{itemize}
It turns out that $T_*$ is not a chain map, although $T^\mu_*$ is.

The first two issues will be handled using a perturbation, as the experts will
have surmised.  The  existence of a chain map
related to $T_*$ is derived from a basic geometric construction, which we
present in the next section.

\subsection{The basic cobordism}

The functoriality of $\HF_*$ with respect to 
oriented cobordisms suggests a method to show that the correspondence $T_*$ is a
chain map.  One  would need (among other things) a
cobordism between $M$ and $\mtau$, over which a representation $\varphi$ and
its cut-and-pasted cousin $\varphi^\tau$ would extend.  We do not know how to
construct such a cobordism, but we can come very close.  We use a variation 
of the idea in \S2 of our earlier paper~\cite{meyerhoff-ruberman:II}
to construct
a cobordism having an additional boundary component, which will be filled in as
the boundary of a $4$-dimensional orbifold.  

Start with a copy of $\Sigma $ sitting inside $M$, so that
$M = A \cup_\Sigma B$.
The labeling of the two components is arbitrary, but having chosen it we
can make the following convention: $\Sigma$ is to be oriented as the boundary
of $A$ (and hence $-\Sigma = \partial B$.)  The mutated manifold $\mtau$ is
then formally a quotient space $(A \coprod B)/\sim$, where `$\sim$'
identifies $x\in \partial B$ with $\tau(x)$ in $\partial A$.    (Since $\tau$
is an involution, it doesn't much matter how we do this, but some care in
making the identifications now will help in the calculations later.)  Fix a
collar neighborhood  $I \times \Sigma $ in $M$.

Consider the manifold $W $ ($= W(M,S^1 \times \Sigma)$ from~\cite{meyerhoff-ruberman:II}) 
constructed as follows:
$$ W = M \times [0,1/4] \cup (I  \times \Sigma) \times  [1/4,3/4]\cup
M \times [3/4,1] 
$$
Cut and paste $W$ along $\Sigma \times [0,1]$, using the involution 
$\tau \times \rm{id}_I$, to obtain a new manifold $\wtau$.\\[2ex]  
\begin{center}
\leavevmode  %%% dr: seems to be necessary for latex to center the figure
\includegraphics{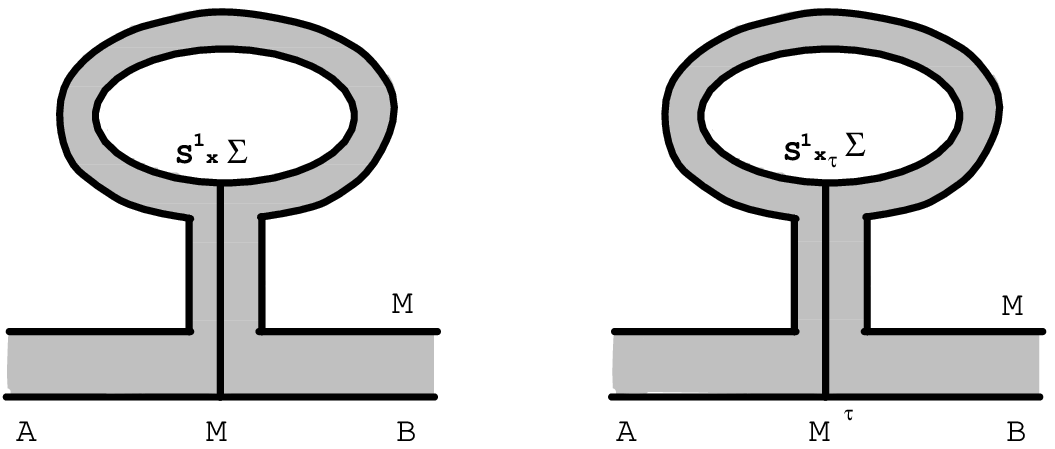}\\
Figure 1
\end{center}

{}From Figure 1, it is clear that the boundary of $\wtau$ consists of a copy of
$M$, a copy of $\mtau$, and a copy of the mapping torus $\mapt$.  From the
point of view of index theory, $\mapt$ will be seen to capture the difference
between $M$ and $\mtau$.  
\begin{lemma}\label{mapt}  The mapping torus $\mapt$ is a Seifert
fibered space $S^2(2,2,2,2,2,2;e)$
over $S^2$, with 6 fibers of multiplicity 2 and Euler class is $-3$.
\end{lemma}
\begin{proof}  Let $\pi :\Sigma \to S^2$ be the branched double cover for
which $\tau$ is the covering transformation.  There are 6 branch points, each
of order 2. Because $\pi\circ \tau = \pi$, the projection extends
to a map $\mapt \to S^2$ which is a fibration away from the fixed points.
The inverse image of each fixed point is a circle, which is covered by nearby
circles with multiplicity $2$.  The Euler class calculation may be done directly,
or more easily by observing that the first homology group of $\mapt$ is
$(\Z/2\Z)^4 \oplus \Z$.  On the other hand, $H_1$ of the Seifert fibered
space $S^2(2,2,2,2,2,2;e)$ is given by $(\Z/2\Z)^4 \oplus \Z/(e+3)\Z$.
\end{proof}

Following the construction${}^{\mathrm{TM}}$ of Fintushel
and Stern~\cite{fs:pseudo}, let $V$ be the mapping cylinder of the projection
$\mapt
\to S^2$, and form an orbifold 
$$
\ytau = \wtau \cup_{\mapt} V
$$
By Lemma~\ref{mapt}, there are $6$ singular points in $\ytau$, each the cone on
$\RP^3$.   It is easily checked that $V$ has the homology of $S^2$, and
that the (rational) self-intersection of the generator of $H_2(V)$ is
trivial. We will also use the notations $V_0$ and $\ytau_0$ for the complement
of an open neighborhood of the six singular points.

We will also need the homology and cohomology groups of $\wtau$.
\begin{lemma}\label{whomology}
The homology and cohomology groups of $\wtau$ are as follows:
$$
H_1 = H^1 = \Z,\quad H_2 = H^2 = \Z^4,\quad H_3 = H^3 = \Z^2.
$$
A choice of basis for $H_1(\Sigma)$ gives a basis for $H_2(\wtau)$. With
respect to a symplectic basis consisting of elements from
$\ker H_1(\Sigma) \to H_1(A)$ or $H_1(B)$, the intersection form is
$$
\begin{pmatrix}
0 & 1\cr 1 & 0\cr
\end{pmatrix}
\oplus\begin{pmatrix}
0 & 1\cr 1 & 0\cr
\end{pmatrix}.
$$
\end{lemma}
\begin{proof}  This all follows from the Mayer-Vietoris sequence, using the
fact that $\tau_* = -1$ on $H_1(\Sigma)$. The calculation of the intersection
form follows~\cite[Proposition 2.1]{meyerhoff-ruberman:II} .
\end{proof}

\subsection{Extending representations}\label{extendsec}

Until further notice, all representations will be in the Lie group G=$\su$ or
$\SO(3)$, so that, for example, $\chi(\pi_1(M))$ refers to the $\su$-character variety
of $M$.  (Since $M$ is a homology sphere, the varieties of $\su$ and $\SO(3)$
representations are the same.)  In this section, we will show that every
$\su$-representation of
$\pi_1(M)$ extends to $\pi_1(\wtau)$ and then to an $\SO(3)$ representation of
$\pi_1^{orb}(\ytau)$.    There always at least two
$\su$  extensions over $\pi_1(\wtau)$, but upon passage to $\SO(3)$ there is only one,
provided that the representation is irreducible when restricted to $\Sigma$.
By restricting to $\pi_1(\mtau)$ we get the one-to-one correspondence between
$\rep(M)$ and $\rep(\mtau)$
referred to above.  It is not obvious (and we will not need to know)
that this is continuous; it is perhaps better to view the maps
$\rep(M) \leftarrow \rep(\wtau) \rightarrow \rep(\mtau)$
as defining  a correspondence in the sense of algebraic geometry.

There is a standard correspondence between $\su$ representations
and flat $\su$ connections, with conjugacy of representations going over to
gauge equivalence.  We will use the two notions interchangeably, without
varying the notation.  In this interpretation, the connection $\alpha^\tau$
is identical to $\alpha$ on $A$, but differs from $\alpha$ on $B$ by a constant
gauge transformation.

Regard $\pi_1(M)$ as being presented as an amalgamation as follows:
$$ \pi_1(M) = \langle \pi_1(A), \pi_1(B) |
i_a(g) = i_b(g) \forall g \in \pi_1(\Sigma) \rangle
$$
Here $i_a,i_b$ are the maps induced by the inclusions of $\Sigma$ into
the two sides.  In this notation, the fundamental group of $\wtau$ is easily
calculated, using van Kampen's theorem.
\begin{equation}\label{pi1w}
\pi_1(\wtau) = \langle \pi_1(A), \pi_1(B), \zz |
\zz^{-1}i_a(g)\zz = i_b(\tau_*(g)) \forall g \in \pi_1(\Sigma) \rangle
\end{equation}
The fundamental group of $\mtau$ has a similar description. 

We will need a more explicit calculation of the effect of $\tau_* $ on $\pi_1(\Sigma)$. 
In figure 2 below, the generators of the fundamental group are $a_1,b_1$ as drawn in the top
half of the surface, and $a_2,b_2$ which are given by $\gamma \alpha_2 \bar{\gamma}$ and
$\gamma \beta_2 \bar{\gamma}$ respectively.  Here, and in what follows, $\bar{x} $ will be used
as a synonym for $x^{-1}$.  
\begin{center}
\includegraphics[scale=.75]{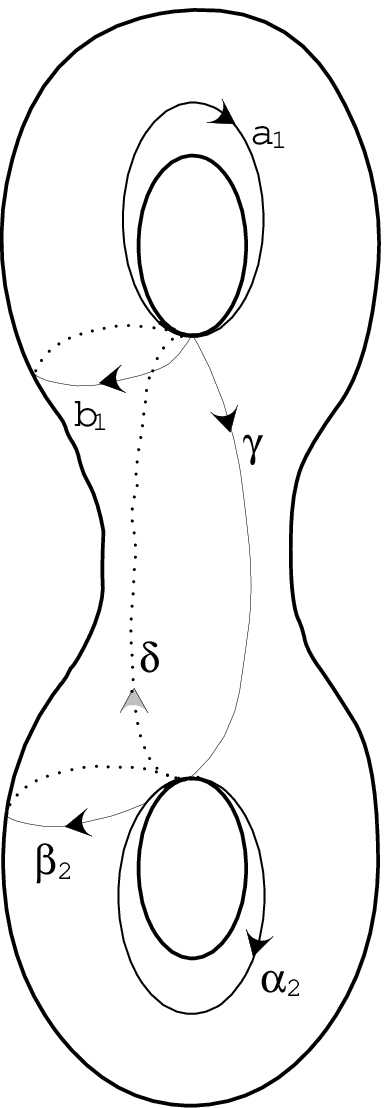}\\
{\bf Figure 2}
\end{center}

Referring to Figure 2, and the curves as labeled therein, we have that
$$
\pi_1(\Sigma) = \langle a_1,b_1,a_2,b_2 |  a_1 b_1 \abar_1 \bbar_1
b_2 a_2  \bbar_2 \abar_2 = 1 \rangle
$$
With respect to these generators, the action of $\tau$ is given by
\begin{align*}
\tau_*(a_1) &= \abar_1\\
\tau_*(b_1) &= a_1 \bbar_1 \abar_1 \\
\tau_*(a_2) &= \bbar_1 b_2 \abar_2 \bbar_2 b_1 \\
\tau_*(b_2) &= \bbar_1 b_2 a_2 \bbar_2 \abar_2 \bbar_2  b_1 
\end{align*} 
For any \su-representation $\varphi$, let $\rm{Stab}(\varphi) = \Z_2, S^1$ or
\su\  be its stabilizer.   
%Note that a representation of $\pi_1(M)$, restricts
%to representations $\alpha_\Sigma, \alpha_A, $ and $\alpha_B$, and that
%$\rm{Stab}(\alpha_A)$ and  $\rm{Stab}(\alpha_B)$ are naturally subgroups
%of  $\rm{Stab}(\alpha_\Sigma)$.
Lemma~\ref{symmrep}, together with the presentation~(\ref{pi1w}) of
$\pi_1(\wtau)$, have the following consequence.

\begin{theorem}\label{extendrep}
Any $\su$ representation $\varphi$ of $\pi_1(M)$ extends to a representation of
$\pi_1(\wtau)$.  The set of  extensions, up to conjugacy, is in one-to-one
correspondence with $\rm{Stab}(\varphi_\Sigma)$.
%$$
%\rm{Stab}(\varphi_\Sigma)
%/\left( \rm{Stab}(\varphi_A)\times \rm{Stab}(\varphi_B) \right)
%.$$
\end{theorem}
\begin{corollary}\label{extendso3rep}
If $\varphi$ is an $\SO(3)$-representation of $\pi_1(M)$ whose restriction to
$\Sigma$ is irreducible, then it has a unique extension to $\pi_1(\wtau)$.
\end{corollary}
\begin{proof}
Since $M$ is a homology sphere, $\varphi$ has a unique lift to an $\su$
representation, which has two extensions to $\pi_1(\wtau)$ according to the
theorem.  But these become the same when projected back to $\SO(3)$.
\end{proof}

When $\mapt$ is filled in to make the orbifold $\ytau$, it is no longer the case
that flat connections extend, because of the possible holonomy around the
$S^1$ fiber.   However, they do extend as flat orbifold connections.
\begin{theorem}\label{extendorb}
Let $\varphi$ be an $\SO(3)$-representation of $\pi_1(M)$ whose restriction to
$\Sigma$ is irreducible. Then it has a unique extension to 
$$\varphi_Y:\pi_1^{orb}(\ytau) =\pi_1(\ytau_0) \rightarrow \SO(3)$$
Furthermore, this representation is non-trivial at each singular point in the
orbifold, and has $w_2(\varphi_Y) $ characterized as follows: it is the unique class in
$H^2(\ytau_0;\Z_2)$ with trivial restriction to $\mapt$ and with non-trivial
restriction to each $\RP^3$ component of $\partial \ytau_0$.
\end{theorem}
\begin{proof}
Because $\tau^2$ is the identity, the presentation~(\ref{pi1w}) implies that
under the extension of $\varphi$ to $\pi_1(\wtau)$, the element $\zz$ goes to an
element of order two in $\SO(3)$.  (The hypothesis on $\varphi_\Sigma$ implies
that $\zz$ can't go to the identity.)  Assuming that the base point was chosen to
be one of the fixed points of $\tau$ on $\Sigma$, the generator of the local
fundamental group near the corresponding singular point of $\ytau$ is $\zz$.  
Hence the representation is non-trivial at that point.  The
generators for the local fundamental groups of the other singular points are all
conjugate to $\zz$, so the representation is non-trivial at each of these points.
That this specifies $w_2(\varphi_Y)$, as described in the statement of the
theorem, may be proved by examining the long exact sequence of the pair
$(V_0,\partial V_0)$.
\end{proof}

In section~\ref{sfsec}, we will compare the Floer-grading of perturbed  flat
connections on $M$ and $\mtau$, using the theory of~\cite{clm:maslov,clm:I,clm:II}.  For
connections which are actually flat, there is a more elementary approach, based on the work
of Atiyah-Patodi-Singer~\cite{aps:I,aps:II} as applied to the operator linearizing the
ASD Yang-Mills equations on $\wtau$.  This approach requires some information about the
\su\ representations of $\pi_1(\mapt)$, and their associated $\rho$ and Chern-Simons
invariants, which we develop in the remainder of this section.   This material may safely
be skipped by those readers who would prefer to pass directly to the proof of the main
theorems on mutation.

\begin{lemma}\label{maptorusreps} The space $\chi(\mapt)$ of \su\ representations of
$\pi_1(\mapt)$ is connected.
\end{lemma}
In the
proof we will consider elements of $\SU(2)$ as unit quaternions.
\begin{proof} The standard presentation of $\pi_1(\mapt)$ as an HNN extension 
is
$$ \pi_1(\mapt) = \langle \pi_1(\Sigma), \zz\, |
\zz^{-1} g \zz  =\tau_*(g)\, \forall g \in \pi_1(\Sigma) \rangle
$$
Thus, a representation $\alpha$ of $\pi_1(\mapt)$ is determined by its restriction
$\alpha_\Sigma$ to $\pi_1(\Sigma)$, and the choice of the element $\zz$ conjugating
$\alpha_\Sigma$ to $\tau^*(\alpha_\Sigma)$.  It is not hard to  show that $\zz$ is
determined up to multiplication by elements of the stabilizer of $\alpha_\Sigma$.  

The restriction map from $\alpha \to \alpha_\Sigma$ is then generically two-to-one, with
connected fibers over the reducible connections on $\Sigma$.  We will show that
the two preimages of an irreducible connection on $\Sigma$ can be  connected via a path in
$\chi(\mapt)$.  The map to 
$\chi^*(\Sigma)$ is in fact a covering space, so the lemma will follow from the fact
that both $\chi^*(\Sigma)$ and the full \su-representation variety $\chi(\Sigma)$ are
connected.  (In fact, $\chi(\Sigma)$  is homeomorphic to $\mathbf{CP}^3$, and the
reducibles are a singular quartic.)

It suffices to find a path connecting  the two preimages (say $\alpha_1,\alpha_{1}'$) of a
single irreducible connection.  The path will be made of 3 pieces: the first connecting
$\alpha_1$ to a reducible connection $\alpha_0$, the second a path in the reducible stratum
{}from  $\alpha_0$ to $\alpha_0'$, and a third connecting
$\alpha_0'$ to
$\alpha_1'$.

Referring to the generators given in Figure 2, let the representation
$\alpha_r, \ r \in [0,1]$ be given by
\begin{align*}
\alpha_r(a_1)& = \exp(-\frac{\pi r}{2}\imath) & 
\alpha_r(a_2) &= \exp(-\frac{\pi r}{2}\imath)\\
\alpha_r(b_1)& = \jmath &
\alpha_r(b_2) &= \jmath\\
\alpha_r(t) &= \exp(\frac{\pi(1- r)}{2}\jmath)
\end{align*}

The representation $\alpha_r$ restricts  to an irreducible representation on $\Sigma$, for
$1 \geq r > 0$.  Note that for these values of $r$, the other representation with the
same restriction to $\Sigma$ differs only in having the opposite sign for $\alpha_r(\zz)$.
On the other hand, 
$\alpha_0$ has
$\alpha_0(\zz) =\imath\jmath=k $, and is reducible on $\pi_1(\Sigma)$, with image lying
in the circle subgroup containing $k $.  Consider the path
$\alpha_{0,s},\, s\in [0,1]$ of representations with the same effect on
$\pi_1(\Sigma)$, but with
$\alpha_{0,s}(\zz) =  \imath(\sin(\pi s) + \cos(\pi s)\jmath).$
Since $ (\sin(\pi s) + \cos\pi s)\jmath)$ is in the centralizer of 
$\alpha_0(\pi_1(\Sigma))$, this gives a path from $\alpha_0$ to $\alpha_0'$, which is the
same as $\alpha_0$ except that $\zz $ is sent to $-k $.  The path $\alpha_r'$, which is
the same as $\alpha_r$ except for reversing the sign of $\alpha_r(\zz)$, provides the third
piece of the path.
\end{proof}

It follows that the variety $\chi(\mapt)$ is has singularities sitting over the
reducibles in $\chi(\Sigma)$, but that the part sitting over the irreducible part of
$\chi(\Sigma) $ is smooth.

Because the Chern-Simons invariant doesn't change on paths of flat connections,
we immediately obtain:
\begin{corollary}\label{maptoruscs}
For any $\alpha \in \chi(\mapt)$, the Chern-Simons invariant $\CS(\alpha) = 0$.
\end{corollary}

The $\rho$-invariant $\rho_{\ad{\alpha_t}}$ is also constant along paths
$\alpha_t$ of flat connections lying in a single stratum of $\chi$, but in
general it may jump as a path descends into a lower
stratum---cf.~\cite{farber-levine:jumps,kirk-klassen:cup}.  The character
variety of $\mapt$ is not smooth, but it turns out that the $\rho$-invariant
doesn't change.  
\begin{lemma}\label{maptorusrho}
For any $\alpha \in\chi(\mapt)$, the invariant $\rho_{\ada}(\mapt) = 0$.
\end{lemma}
\begin{proof}
Because the $\rho$-invariant is locally constant on strata, it suffices to check the
vanishing for a single representation in each component of the three strata of
$\chi(\mapt)$.   The technique is the same in each stratum, so we just check
the case when the restriction to $\mapt$ is irreducible.  Let $\alpha$ be a
representation sending
$a_1 \to \imath$, $a_2 \to \jmath$, and $b_1, b_2 \to 1$.  
There are two choices $(\pm k)$ for $\alpha(\zz)$; the argument works with either
one. Notice that  $\alpha_\Sigma$ extends over the obvious genus-2 handlebody
$C$ with boundary
$\Sigma$, so that $\alpha$ extends over the ($4$-dimensional) mapping torus 
$ S^1\times_\tau C$.  

Essentially by definition~\cite{aps:II} 
$$ \rho_{\ada}(\mapt) =  3\,\text{sign}(S^1\times_\tau C) -
\text{sign}(S^1\times_\tau C;\ada) $$ Since $\tau_* = -1$ on $H_1(C)$, a simple
Wang sequence shows that $H_2(S^1\times_\tau C) = 0$, so the first
signature vanishes.  To compute the second signature, we compare the cohomology
of
$S^1\times_\tau C$ with that of $\mapt$.  Both of these 
are computed via Wang sequences, summarized in the
following diagram in which  $\ada$-coefficients are understood.
$$
\begin{CD}
0 @>>> H^1(\mapt)  @>>>  H^1(\Sigma) 
@>\tau^* -1>> H^1(\Sigma)   @>>>    H^2(\mapt)   @>>>    0\\
@NNN @AAA @AAA @AAA @AAA @NNN\\
0 @>>> H^1(S^1\times_\tau C)  @>>>  H^1(C) 
@>\tau^* -1>> H^1(C)   @>>>   H^2(S^1\times_\tau C)  @>>>    0
\end{CD}
$$
Recalling that $\alpha$ is irreducible on $\Sigma$, we have that
$H^2(\Sigma;\ada) = H^0(\Sigma;\ada)$ vanishes.  Similarly, $H^2(C;\ada) =
H^0(C;\ada)$.  One can compute that $\tau^*$ is the identity, using group
cohomology.  Alternatively, since $H^1(\Sigma;\ada) $ is the tangent space to
$\chi(\pi_1(\Sigma))$, on which $\tau$ acts by the identity,  $\tau^* = id$.
A similar remark applies to the action on $H^1(C;\ada) $. Since $C$ has
the homotopy type of $1$-complex, 
$0 = H_2(C;\ada) \cong H^1(C,\Sigma;\ada) $,
and  $H^1(C;\ada) \to H^1(\Sigma;\ada)$ is an injection.
A diagram chase shows that
$H^2(S^1\times_\tau C;\ada)$ injects into $H^2(\mapt;\ada)$, which implies that
the second signature vanishes as well.
\end{proof}

The vanishing of $\CS(\alpha)$ and $\rho_{\ad\alpha}$ could equally well
have been obtained using the Seifert-fibered structure of $\mapt$, as
in~\cite{fs:instanton}.  Alternatively, the jumps in the $\rho$-invariant
could presumably be calculated using the techniques
of~\cite{farber-levine:jumps,kirk-klassen:cup}.

\subsection{Perturbations}\label{pertsec}

It is not necessarily the case that the flat connections on $M$ form a smooth
$0$-dimensional variety.  If they don't, then it is necessary to perturb the
flatness equations, in order to define the Floer homology groups. For the
purposes of this paper, it is desirable to make the perturbations in such a
way that perturbed-flat connections can be cut and pasted along $\Sigma$.  Now
a standard method (compare~\cite{taubes:casson,floer:instanton}) for perturbing
the equation
$F_A= 0$ is to replace the right side by an $su_2$-valued
$2$-form supported in a neighborhood of a link
$L \subset M$.  A suitable class of $2$-forms can be defined as follows:
Choose first a link $L$ in $M$, and for each component $L_i$ a $C^2$ function
$\bar{h}_i: [-2,2] \to \R$.  The whole collection determines a function on
the space of connections
$$
h(A) = \sum_i\int_{D^2}\bar h_i (\mathrm{tr}\,\mathrm{hol}_{L_i}(x,A)) \eta\,dx
$$
where $\eta$ is a bump function on the normal disk $D^2$ and
$\mathrm{hol}_{L_i}(x,A)$ denotes the holonomy around a curve in $L_i \times
D^2$.

Any such $h$, as a function on the space of connections, has a gradient
$\nabla h(A)$, which is naturally a 1-form on $M$, supported near $L$.
The solutions to the equation $F_A = *\nabla h(A)$ are called \hflat connections, and it
can be shown that for a sufficiently complicated link
$L$, a generic choice of
$h$ will result in a smooth moduli space $\chi_h$ of \hflat connections.   By construction,
the \hflat connections are the critical points of the function $\CS_h = \CS - h$.

In order to cut/paste an \hflat connection along $\Sigma$, it is
necessary that the link $L$ along which the perturbation is supported be
disjoint from $\Sigma$. (In other words, if $L$ hits $\Sigma$, the \hflat connections on
$\Sigma$ don't enjoy a symmetry property analogous to~\ref{symmrep}.)
Recall that $M$ is divided into two pieces $A$ and $B$.  If $L$ is disjoint
{}from $\Sigma$, then we can consider separately the \hflat connections on $A$
and $B$.  Such connections can clearly be glued up exactly when they agree on
$\Sigma$.  In other words,  $\chi_h(M)$ is (essentially) the fiber
product of 
$$\chi_h(A) \times_{\chi(\Sigma)} \chi_h(B) $$
(This description must be modified, in a standard way, when connections are
reducible along $\Sigma$.)  In particular, the correspondence between flat
connections on $M$ and $\mtau$ continues to 
hold for
\hflat connections which are flat along
$\Sigma$.

\begin{theorem}\label{perturb}
Let $M$ be a homology sphere, and let $\Sigma$ be an embedded genus-$2$ surface. 
Then there is a link $L$ in the complement of $\Sigma$, and a
perturbation of the equation
$F_A = 0$ to an equation of the form $F_A  =*\nabla h(A)$, where the
2-form $*\nabla h(A)$ is supported near $L$, with the following properties:  For
any solution $\alpha$ (i.e.~\hflat connection), 
$\HH^*(M;\ada)= 0$, and the restriction of $\alpha$ to $\Sigma$ is an irreducible
flat connection.
\end{theorem}
Here $\HH^j(M;\ada)$ denotes the space of $\ada$-valued harmonic $j$-forms
on $M$.

We now turn to setting up the proof of Theorem~\ref{perturb}.
Let $r_A$ (resp. $r_B$) denote the restriction of connections from $A$
(resp. $B$) to $\Sigma$. 
Finding a perturbation $h$ for which
$\chi_h(M)$ is  smooth breaks into two steps: smoothness for the two
sides, and transversality of the maps $r_A$ and $r_B$.  We will treat these
issues using results from the paper of  C.~Herald~\cite{herald:perturbations}.

Following that paper, let $\M(X)$ denote the flat $\SU$-connections on a manifold 
$X$, and $\MZ$ (resp. $\MU,\, \MSU$) denote the connections with stabilizer
$\Z_2$ (resp. $\mathrm{U}(1),\, \SU(2)$).  If $X$ is a manifold with boundary,
and
$G \subset H \subset \SU(2)$, then $\M^{G,H}$ consists of connections
with stabilizer (on $X$) equal to  $G$ and stabilizer of the restriction to
$\partial X$ equal to $H$.   For a generic perturbation $h$, $\M_h(X)$ will
be a stratified space, with strata indexed by the various possible pairs
$(G,H)$.  In the case when $X= A$ or $B$, so that $\partial X =
\Sigma$ is a surface of genus $2$, not all possible combinations of $(G,H)$
occur as strata of
$\M_h(X)$, after the
perturbation: Only the  $(\Z_2,\Z_2),\ (\Z_2,U(1)),\ 
(U(1),U(1)),\ \mathrm{and}\ (\SU(2),\SU(2))$ strata will appear.  

The paper of Herald gives some additional results pertaining to the
restriction map $r_A$ from the various strata $\M_h^{G,H}(A)$ to
$\M_h^{H}(\Sigma)$.  In rough terms,
Herald's paper shows that, for
generic  perturbations of the equations on $A$ and
$B$, the maps $r_{A}$ and $r_{B}$ are transverse.  Of course, since all the
spaces involved are, at best, stratified, the transversality must be taken in a
suitable sense.  Precise statements are given in Lemma 33 and Proposition 34
of~\cite{herald:perturbations}.  From these results, we will prove:
\begin{lemma}\label{strtrans}
Suppose that $M$ is a homology sphere.  Then for a generic perturbation $h$ of
the equations, the images of the restriction maps 
$$
\begin{CD}
\M_h^{G,H}(A) @>r_A>> \M^H(\Sigma) @<r_B<< \M_h^{G',H}(B)
\end{CD}
$$
are empty, except if $(G,H) = (G',H) = (\Z_2,\Z_2)$ or $(\SU(2),\SU(2))$.
In the $(\Z_2,\Z_2)$ case, the maps will be transverse, while in the 
$(\SU(2),\SU(2))$ case the intersection is isolated at the trivial connection.
\end{lemma}
\begin{proof}
According to~\cite[Theorem 15]{herald:perturbations}, generic perturbations of
the equations on $A$ and $B$ will eliminate all but the $(\Z_2,\Z_2),\ 
(\Z_2,\U(1)),\ (\U(1),\U(1)), $ and $(\SU(2),\SU(2))$ strata.  Evidently,
connections with the same restriction to $\Sigma$ must have the same stabilizer
on $\Sigma$, so there are a limited number of cases to examine.  Here is the
full list of possibilities (modulo switching the letters $A$ and $B$), together
with a description of what happens in each case:\\[2ex]
\begin{enumerate}
\item  $r_A(\M_h^{\SU(2),\SU(2)} (A)) \cap r_B(\M_h^{\SU(2),\SU(2)} (B)) =
\{\Theta$\}:\hspace{1em}
When $h=0$, the only points in the fiber product would be the flat
connections with image of the holonomy in $\Z_2$.  But since $M$ is a homology
sphere, the only possibility is the trivial connection, which is an isolated
point in
$\chi(M)$ (cf.~\cite{akbulut-mccarthy}).  This situation is stable under small
perturbations, so remains true for the perturbed moduli space.\\[1.5ex]
\item $r_A(\M_h^{\U(1),\U(1)} (A)) \cap r_B(\M_h^{\U(1),\U(1)} (B))=
\emptyset$:\hspace{1em}   This intersection is empty before the
perturbation, since $M$ is a homology sphere, and so remains empty if $h$ is
sufficiently small.\\[1.5ex]
\item $r_A(\M_h^{\Z_2,\U(1)} (A)) \cap r_B(\M_h^{U(1),U(1)} (B))
= \emptyset$:\hspace{1em}
The image of $\M_h^{U(1),U(1)} (B)$ is a 2-dimensional submanifold of  the
smooth, 4-dimensional manifold 
$\M^{\U(1)}(\Sigma)$.
Theorem 15 of Herald's paper says that if h is generic, then $\M^{\Z_2,\U(1)}
(A)$ is
$0$-dimensional, so its image under $r_A$ is a finite set of points in 
$\M^{\U(1)}(\Sigma)$.
 Moreover, Lemma 33 of~\cite{herald:perturbations} says that the restriction
map from an appropriate subset of $\{\mathrm{connections}\} \times
\{\mathrm{perturbations}\}$ to $\M^{\U(1)}(\Sigma)$ is a submersion.   This
means that $r_A(\M_h^{\Z_2,\U(1)} (A))$ can be moved
arbitrarily in 
$\M^{\U(1)}(\Sigma)$, and so can be arranged to miss $r_B(\M_h^{\U(1),\U(1)}
(B))$.\\[1.5ex]
\item  $r_A(\M^{\Z_2,\U(1)} (A)) \cap r_B(\M^{\Z_2,U(1)} (B))= \emptyset$:\hspace{1em}
This is the same argument: the two  finite sets of points can be made disjoint
by a small perturbation.\\[1.5ex]
\item $r_A(\M_h^{\Z_2,\Z_2} (A)) \trans r_B(\M_h^{\Z_2,\Z_2} (B)) = $ a finite
set of points:\hspace{1em}
As in the previous two items; the transversality follows from~\cite[Lemma
33]{herald:perturbations}.
\end{enumerate}
\end{proof}

\begin{proof}[Proof of Theorem~\protect{\ref{perturb}}]
Choose a perturbation $h$ as described in
Lemma~\ref{strtrans}.  Then $\M_h(M)$, being the fiber product of $\M_h(A)$ and
$\M_h(B)$, will consist of the (isolated) trivial connection, together with the
fiber product of the $(\Z_2,\Z_2)$ strata.  A final application of~\cite[Lemma
33]{herald:perturbations} shows that the maps $r_A$ and $r_B$ are generically
immersions of smooth $3$-manifolds into into $\MZ(\Sigma)$, so the fiber
product is a finite set of points.  The maps on harmonic forms induced by $r_A$
and
$r_B$ are just the differentials of those maps.  So the fact that the cohomology
at any point in the fiber product is trivial is simply a restatement of the
transversality condition.
\end{proof}

As described in~\cite[\S3.3]{braam-donaldson:knots}, the perturbed
equations on $M$ extend to perturbed anti-self-duality equations on
$\ytau$.  Here, the perturbation is supported on a neighborhood of $L \times \R
\subset \ytau_\infty$, and is hence a `time independent' deformation of the
ASD equations, in the terminology of~\cite{braam-donaldson:knots}.  (The notation
$\ytau_\infty$ indicates that half-infinite tubes are added
along the boundary components of $\ytau$.)  In our situation, the solutions will be flat on
the complement of the neighborhood
$\nu(L \times \R)$, and so we will continue to refer to them as \hflat connections.   
\begin{defn}
The moduli space of $h$-ASD orbifold $\SO(3)$ connections on $\ytau_\infty$, with
exponential decay to  \hflat connections $\alpha,\beta$ on $M$ and $\mtau$, respectively,
will be denoted by
$$
\M_h(\ytau;\alpha,\beta)
$$
If, as in the corollary below, the
Stiefel-Whitney class $w$ or the charge
$c$ is specified in advance, the notation will be expanded to $\M_{h,c}^w$.
\end{defn}
\begin{corollary}\label{pertmutation} Let $M$ be a homology sphere containing a
genus-$2$ surface $\Sigma$.  Let $ *\nabla h(A) $ be a 2-form supported
along a link $L$ in the complement of $\Sigma$, having the properties described
in Theorem~\ref{perturb}.  Then for any \hflat connection on $M$,
there is a unique extension to an \hflat orbifold connection on $\ytau.$ 
(The behavior at the singular points, and $w= w_2$ of the bundle are as
specified in Theorem~\ref{extendorb}.)   Restricting to the boundary gives
a one-to-one correspondence 
$$
\M_h(M) \leftrightarrow \M_{h,0}^w(\ytau) \leftrightarrow\M_h(\mtau).
$$
\end{corollary}
\begin{proof}[Proof of corollary]
Apply the constructions of section~\ref{extendsec} to the $3$-manifold $N$
obtained by removing a tubular neighborhood of $L$ from $M$.  (Note that the
arguments in section~\ref{extendsec} did not assume that $M$ was closed, and
hence apply without change to $N$.)  The result is a cobordism from $N$ to
$N^\tau$, obtained by removing a tubular neighborhood of $L \times \R$ from
$\ytau$.   By Theorem~\ref{extendorb}, there are 
one-to-one correspondences between $\chi(N)$, the flat orbifold connections on
$\ytau - \nu(L \times \R)$, and $\chi(N^\tau)$.  The argument that a flat
connection on $M - \nu(L)$ extends uniquely to an \hflat connection on $M$
applies in the $4$-dimensional situation, and shows that a flat connection on
$\ytau - \nu(L \times \R)$ extends uniquely to an \hflat connection on $\ytau$.
\end{proof}

The \hflat extension defines the unique element of $ \M_{h,0}^w(\ytau;\alpha,\alpha^\tau)$.
We will show later that the formal dimension of this moduli space is $0$; we need to
know that the one point in the moduli space is a smooth point.  The usual deformation
theory says that this will follow from the following lemma.
\begin{lemma}\label{smooth} 
Let $A$ be the unique element of  $\M_{h,0}^w(\ytau;\alpha,\alpha^\tau)$.  Then the
space of harmonic forms
 $\HH^1(\ytau;\ad A)$ vanishes.
\end{lemma}
\begin{proof}
In the unperturbed situation (i.e. if $\alpha$ were a smooth isolated flat connection), this
could be readily proved by interpreting 
$\HH^1(\ytau_\infty;\ad A)$ as a twisted cohomology group, and then computing topologically,
using the isomorphisms $H^1(M;\ada) \xleftarrow{i^*} \HH^1(\ytau;\ad A) \xrightarrow{i^*} 
H^1(\mtau;\ada^\tau) = 0$.  Equivalently, one could interpret the cohomology group in terms
of group cohomology, and use the presentation~(\ref{pi1w}) of
$\pi_1(\wtau)$ to obtain the same isomorphism.  Note that either of these arguments would
apply to  $\ytau - \nu(L \times \R)$.  The proof of the lemma would thus be completed if
there were an appropriate Mayer-Vietoris principle for $\ad(A)$-valued harmonic forms on
manifolds such as $\ytau_\infty$.  Here it would have to be applied to the decomposition of
$\ytau$ into $\ytau - \nu(L \times \R) \cup \nu(L \times \R)$.  While such an argument would
undoubtedly succeed, we know of no convenient reference, and proceed somewhat differently.

It is explained carefully in Appendix A of~\cite{clm:I} that the (exponentially decaying)
harmonic forms on $\ytau_\infty$ can be identified with harmonic forms on $\ytau$ (regarded
as having a long, but finite, collar on its boundary.)   Double $\ytau$ along its boundary,
to obtain an orbifold $Z$, which contains a copy of $\nu(L \times S^1)$.  By a
Mayer-Vietoris argument, found in~\cite[Appendix B]{clm:I}, the harmonic forms
$ \HH^1(\ytau;\ad A) \cong \HH^1(Z;\ad A)$, so it suffices to show the latter.  But now we
can use the Mayer-Vietoris principle for
$$
Z = (Z - \nu(L \times S^1)) \cup \nu(L \times S^1)
$$
As remarked above, the restriction
$\HH^1(Z - \nu(L \times S^1);\ad A) \rightarrow \HH^1(M -\nu(L);\ad A)$ is an isomorphism,
as is the restriction 
$\HH^1(\nu(L \times S^1);\ad A) \rightarrow \HH^1(\nu(L);\ad A)$.
By a diagram chase, the vanishing of the harmonic forms for $Z$ follows from the same on
$M$, plus the Mayer-Vietoris principle for forms on $M =(M - \nu(L)) \cup \nu(L)$.
\end{proof}

\subsection{Spectral Flow}\label{sfsec}
In this section we show that corresponding \hflat connections in $\chi_h(M)$
and $\chi_h(\mtau)$ have the same grading in Floer homology.  As remarked at the end of
section~\ref{extendsec}, for connections which are actually flat, this can be proved using
the Atiyah-Patodi-Singer index theorem.  (One would need, in addition to the $\rho$ and
$\CS$ invariants computed there, a calculation of the twisted signature of $\wtau$.)
This approach is not workable for the
\hflat connections because it would involve a direct calculation of the index of a
Dirac-type operator, whose kernel and cokernel are not topological invariants. 
An alternative technique, which we adopt, is to use the definition of the grading
in terms of spectral flow of paths of differential operators.    Given the
splitting of $M$ along
$\Sigma$ into pieces
$A$ and $B$, work of Cappell, Lee, and Miller~\cite{clm:maslov,clm:I,clm:II} (and
others~\cite{daniel:splitting,nicolaescu:sf}) calculates the spectral flow on
$M$ in terms of  spectral flows of operators on $A$ and $B$, and an `interaction
term' called the Maslov index.  We will show that the these terms do not
change in passing from $\alpha$ to $\alpha^\tau$.   

Suppose that $M$ has a Riemannian metric, which for future reference will be
chosen to be a product $\Sigma \times [-r,r]$ near $\Sigma$.  Here $r$ is
chosen sufficiently large, so that Theorem C of~\cite{clm:II} applies to all
paths of connections under consideration.  Let
$\alpha$ be a connection on
$M$ and choose a path
$\alpha_t$ of connections from the trivial connection $\Theta_M$ to $\alpha$. 
For each $\alpha_t$ there is defined a first-order elliptic operator 
\begin{equation}
D_{\alpha_t} = 
\begin{pmatrix}
0 & d_{\alpha_t}^*\\
d_{\alpha_t} & -*d_{\alpha_t}
\end{pmatrix}:
\begin{matrix}
\Omega^0(M;\ada_t) \\ 
\oplus\\
 \Omega^1(M;\ada_t)
\end{matrix}
\xrightarrow{\quad}
\begin{matrix}
\Omega^0(M;\ada_t)\\
 \oplus\\
 \Omega^1(M;\ada_t)
\end{matrix}
\end{equation}
The Floer grading of $\alpha$ is then given by the spectral flow of the path of
operators $\{D_{\alpha_t}| t \in [0,1]\}$.  To be more precise, it is the
$(-\epsilon,\epsilon)$--spectral flow in the terminology of~\cite{clm:II},
i.e.~the intersection number of the graph of the eigenvalues of $D_{\alpha_t}$
in $[0,1] \times \R$ with the line from $(0,-\epsilon)$ to $(1,\epsilon)$.   The
work of Cappell--Lee--Miller actually computes the
$(\epsilon,\epsilon)$--spectral flow, which we will denote by $\SF^\epsilon$. 
Fortunately, it is not hard to account for the difference between the two
spectral flows.

\begin{lemma}\label{epsilon} Let $\epsilon >0$, and suppose that $D_t$ be a path
of self-adjoint operators, such that any nonzero eigenvalue $\lambda$ of $D_0$
or $D_1$ satisfies $\lambda > \epsilon$.  Then the
$(-\epsilon,\epsilon)$ spectral flow $\SF(D_t)$ and the $(\epsilon,\epsilon)$
spectral flow $\SF^\epsilon(D_t)$ are related by
$$
\SF(D_t) = \SF^\epsilon (D_t) + \dim(\ker(D_0))
$$
\end{lemma}
\begin{remark}
As
in~\cite{clm:II}, the hypothesis of the lemma will hold if $\epsilon $ is chosen
to be $1/r^2$ for a sufficiently large value of $r$.
\end{remark}
\begin{proof}
The hypothesis implies that a $0$-eigenvalue of $D_0$ either flows to a
positive eigenvalue $> \epsilon$ of $D_1$, to a $0$-eigenvalue, or to a
negative eigenvalue $< -\epsilon$.  In the first case, it contributes $+1$ to
$\SF^\epsilon$ and $0$ to $\SF$.  In the latter two cases, it contributes $-1$
to $\SF$, and $0$ to $\SF^\epsilon$, and the result follows since all
other paths of eigenvalues contribute the same to both counts of spectral
flow.  Another way to phrase the argument is that any path must have $0$
intersection number of with the (oriented) triangle in the $(t,\lambda)$-plane
with vertices
$(0,-\epsilon)$, $(0,\epsilon)$, and $(1,\epsilon)$.  From this point of view,
the lemma states the obvious fact that the intersection with the piece along
the $\lambda$-axis is given by $\dim(\ker(D_0))$.
\end{proof}

The spectral flow $\SF_M(D_\Theta,D_\alpha)$ being independent of the choice of
path means that we can choose a path of connections which is well-suited to
cutting and pasting.  
\begin{construction}\label{path}
Recall from Theorem~\ref{perturb} that the \hflat connection $\alpha$ is flat and
irreducible on $\Sigma$.  Since
$\chi(\Sigma)$ is connected, there is a smooth path $\alpha_{\Sigma,t}$ of flat connections on
$\Sigma$, with $\alpha_{\Sigma,0}$ the trivial connection.   It is convenient to choose
the path so that $\alpha_{\Sigma,t}$ is irreducible for
$0 < t \leq 1$.  This may be done since the space of irreducible flat
connections on $\Sigma$ is connected.  One consequence of this choice is that
the kernel of the `tangential' operator
$\hat{D}_{\alpha_t}$, given by
$H^0(\Sigma,\ada_t) \oplus H^0(\Sigma,\ada_t) \oplus H^1(\Sigma,\ada_t)$, is constant
for $t > 0$.  Using a partition of unity, extend $\{\alpha_{\Sigma,t}\}$ to a path
of connections on $M$ with $\alpha_0$ the trivial connection.   The connections may
be assumed to be flat, and pulled back from $\Sigma$, on the tube $[-r,r]
\times \Sigma$.

By Lemma~\ref{symmrep}, there is an element $g_t \in \su$ conjugating
$\alpha_{\Sigma,t}\circ\tau_*$ to $\alpha_{\Sigma,t}$.  Since $\alpha_{\Sigma,t}$ is
irreducible, $g_t$ is determined up to sign, and so the induced path in $\SO(3)$ is smooth. 
By path-lifting for the double covering $\su \rightarrow \SO(3)$,  the path
$g_t$ may thus be assumed to be smooth.  By choosing the path to agree with some
previously specified model path near $t=0$, we can arrange that the path has a smooth
extension to $t=0$.
\end{construction}

The result of this construction is that each element in the path $\alpha_t$ can be cut
and pasted to give a smooth path of connections $\alpha_t^\tau$ on $\mtau$, giving
rise to paths of operators  $D_{\alpha_t}$ and $D_{\alpha_t^\tau}$ on $M$ and
$\mtau$. The splitting
technique for  computing
the spectral flow of $D_{\alpha_t}$ and $D_{\alpha^\tau_t}$ involves the restriction
of $su_2$-valued forms (and operators on these spaces) on $M$ to those on $A$ and $B$
and subsequently to $\Sigma$. 
We make the following convention:
\begin{convention}\label{ident}
Forms on $A$, $B$
and $\Sigma$, all viewed as submanifolds of $\mtau$, are identified with the same forms
when $A$, $B$ and $\Sigma$ are viewed as submanifolds of $M$.
\end{convention}
The difference between gluing
$A$ to $B$ via $\tau$ and via the identity is encoded in the restriction maps from 
forms on $A$ or $B$ to those on $\Sigma$.   By construction,
$A$ as a submanifold of $\mtau$ is identified with $A$ as a submanifold of $M$,  and
under this identification, $\alpha^\tau_A $ is identified with $\alpha_A$.  In this
way, the restriction map
$$ \Omega^*(A;\ada^\tau_t) \rightarrow \Omega^*(\Sigma;\ada^\tau_t)  $$
is the same as the restriction map 
$$ \Omega^*(A;\ada_t) \rightarrow \Omega^*(\Sigma;\ada_t).  $$
In contrast, the
restriction map from forms on $B$ (viewed as a subset of $\mtau$) to forms on
$\Sigma$ is given by the composition
\begin{equation}
\Omega^*(B;\ada_t) \xrightarrow{r_B^*} \Omega^*(\Sigma;\ada_t) \xrightarrow{\tau^*}
\Omega^*(\Sigma;\ada_t) 
\end{equation}
where $r_B^*$ is the restriction map from forms on $B \subset M$  to forms on
$\Sigma$ and $\hat\tau^*$ is induced by $\tau$ as described in equation~(\ref{taut}) at
the beginning of section 1. (Note that $\hat\tau^*$ actually depends on $t$, but this
dependence will be suppressed in the notation.)

The operators
$D_{\alpha_t}(A)$ and $D_{\alpha_t^\tau}(A)$, obtained by restricting from $M$ and
$\mtau$ to $A$ are identical, as are the the restrictions 
$D_{\alpha_t}(B)$ and  $D_{\alpha_t^\tau}(B).$    It is important to remark, however, 
that because of the action of $\hat\tau^*$ on $\Omega^*(\Sigma;\ada)$, the
boundary-value problems on $B$ associated to the two operators are not {\sl a priori}
the same.

We now summarize the splitting results from~\cite{clm:I,clm:II} which will be used
to compare the spectral flow of the path
$\alpha_t$ with that of its mutated cousin.  As described in the paragraphs before
Theorem C in ~\cite{clm:II}, divide the interval $[0,1]$ into sub-intervals $0 = a_0 <
a_1
\ldots a_n = 1$ with the property that for $ t \in [a_i,a_{i+1}]$ there are no
eigenvalues of
$\hat{D}_{\alpha_t}$ in the intervals $(K_i,K_i+\delta)$ and $(-K_i-\delta,-K_i)$,
for some positive
$K_i$ and $\delta$.   The spectral flow $\SF^\epsilon$ is then the sum of the spectral
flows on the subintervals, so it suffices to compare on the intervals $[a_i,a_{i+1}]$.

For $t$ in such an interval, there are smoothly varying 
Atiyah-Patodi-Singer type boundary conditions for the $D_{\alpha_t}(A)$ and
$D_{\alpha_t}(B)$  (i.e.~the restriction of $D_{\alpha_t}$ to $A$
(resp.~$B$)), described as follows.  The finite dimensional space
$\mathcal{H}(t,K_i) \subset \Omega^0(\Sigma;\ada_t) \oplus \Omega^1(\Sigma;\ada_t)$
is defined to be the span of the eigenfunctions of $\hat D_{\alpha_t}$ with eigenvalue less
than $K_i$ in absolute value.  Note that the images of $\ker(D_{\alpha_t}(A))$
(resp.~$\ker(D_{\alpha_t}(B))$) under restriction to $\Sigma$ give Lagrangian
subspaces $L_t(A)$ (resp.~$L_t(B)$) in
$\ker(\hat D_{\alpha_t})$.  These are extended to Lagrangian subspaces
of $\mathcal{H}(t,K_i)$ by defining
\begin{align*}
\mathcal{L}_t(A) & = L_t(A) \oplus \left[P_+(t) \cap \mathcal{H}(t,K_i) \right]\\
\mathcal{L}_t(B) & = L_t(B) \oplus \left[P_-(t) \cap \mathcal{H}(t,K_i) \right]
\end{align*}
where $P_{\pm}(t)$ are the sums of the positive/negative eigenspaces of $\hat D_t$.
The resulting path (for $t\in [a_i,a_{i+1}]$) of pairs of Lagrangian
subspaces defines (cf.~\cite{clm:maslov}) a Maslov index
$$\mathrm{Mas}(\mathcal{L}_t(A),\mathcal{L}_t(B)  ) \in \Z.$$

The operator $D_{\alpha_t}(A)$
(for $t \in [a_i,a_{i+1}]$) is Fredholm on the domain consisting of forms in 
$L_1^2(\Omega^0(A;\ada_t) \oplus \Omega^1(A;\ada_t))$ whose restriction to $\Sigma$
lie in 
$$
\mathcal{L}_t(A) \oplus P_+(t,K_i).
$$
Here $ P_+(t,K_i)$ denotes the span of the eigenfunctions of $\hat D_{\alpha_t}$ with
eigenvalue greater than  $K_i$.  Similarly, $D_{\alpha_t}(B)$ is Fredholm, where now its
domain is specified by the requirement that the forms, upon restriction to $\Sigma$,
lie in $ \mathcal{L}_t(B) \oplus P_-(t,K_i).$  It follows that the spectral flows
$\SF^\epsilon(D_{\alpha_t}(A))$ and $\SF^\epsilon(D_{\alpha_t}(B))$ are defined for $t\in
[a_i,a_{i+1}]$.  

\begin{theorem}\label{sfequal}
For any irreducible \hflat connection $\alpha$, $\SF_M(D_\Theta,D_\alpha) =
\SF_{M^\tau}(D_\Theta,D_{\alpha^\tau})$.  Equivalently,  the Floer grading of
$\alpha^\tau$ is equal to the Floer grading of $\alpha$. 
\end{theorem}
\begin{proof}
Because $M$ and $\mtau$ are both homology spheres, and $\alpha_0$ and $\alpha_0^\tau$
are trivial connections, Lemma~\ref{epsilon} implies that it suffices to show that
 $\SF^\epsilon_M(D_\Theta,D_\alpha) =
\SF^\epsilon_{M^\tau}(D_\Theta,D_{\alpha^\tau})$.  Moreover, the discussion in the
preceding paragraphs implies that it suffices to prove this equality when $t$ ranges
over the interval $[a_i,a_{i+1}]$.

Theorem C of~\cite{clm:II} states that for sufficiently small $\epsilon$,
\begin{multline*}\tag{$*$}\label{clmsf}
\SF^\epsilon_M(D_{\alpha_t}) =
\SF^\epsilon(D_{\alpha_t}(A)) + \SF^\epsilon(D_{\alpha_t}(B))\\ +
\mathrm{Mas}(\mathcal{L}_t(A),\mathcal{L}_t(B)) 
+\frac{1}{2}[\dim \ker \hat D(a_{i+1}) - \dim \ker \hat D(a_{i})]
\end{multline*}
where all of the terms are calculated on the interval $[a_i,a_{i+1}]$.
Hence it suffices to show that the terms on the right hand side in the analogous
formula $(*^\tau)$ for $\SF^\epsilon_{M^\tau}(D_{\alpha_t^\tau})$ are the same as
those above.  

Using the convention~\ref{ident}, the operator $\hat D_{\alpha^\tau}$, when
$\Sigma$ is viewed as a submanifold of $\mtau$ is identified with $\hat D_\alpha$
(for $\Sigma$ viewed as a submanifold of $M$).  Hence the kernel of $\hat D_\alpha$
is unchanged when $\alpha $ is replaced by $\alpha^\tau$, so the last terms  in
equations~\eqref{clmsf} and $(*^\tau)$ are the same.  Similarly,  the Lagrangian
subspace $\mathcal{L}_t(A)$ does not change, whether $A$ is viewed as a submanifold
of $M$ or of $\mtau$.  It follows that $\SF^\epsilon(D_{\alpha_t}(A)) =
\SF^\epsilon(D_{\alpha^\tau_t}(A)) $, because the two refer to spectral flow of
operators which are viewed as identical. 

To show that $\SF^\epsilon(D_{\alpha_t}(B))$ and the Maslov index term do not
change under the mutation, it suffices to show that $\hat\tau^*$ takes
$\mathcal{L}_t(B)$ to itself.  Now $\hat\tau^*$, being an automorphism of
$\alpha_t$, commutes with
$\hat D_{\alpha_t}$, and hence preserves the eigenspaces of that operator.  In
particular, the summand $\left[P_-(t) \cap \mathcal{H}(t,K_i) \right]$ of
$\mathcal{L}_t(B)$ is preserved by $\hat\tau^*$, so we need to know the effect of
$\hat\tau^*$ on the Lagrangian subspace $L_t(B) \subset \ker (\hat D_{\alpha_t})$.
\begin{claim}%\label{tau*}
Let $\alpha$ be a flat connection on $\Sigma$.  If $\alpha$ is irreducible, so
that $\ker(\hat{D}_{\alpha}) \cong H^1(\Sigma;\ada)$, then $\tau^*$ acts as the
identity on $\ker(\hat{D}_{\alpha})$.  If $\alpha$ is $\su$-reducible, then
$$
\ker(\hat{D}_{\alpha}) \cong (H^0(\Sigma) \otimes su_2) \oplus ( H^0(\Sigma)
\otimes su_2) \oplus (H^1(\Sigma) \otimes su_2)
$$
and $\tau^*$ is the identity on the first two summands and $-1$ on the third.
\end{claim}
\begin{proof}[Proof of Claim]
\renewcommand{\qedsymbol}{{$\openbox_{claim}$}}
In the case that $\alpha_\Sigma$ is irreducible, one can show that
that $\hat\tau^*$ is the identity, using group
cohomology.  Alternatively, since $H^1(\Sigma;\ada) $ is the tangent space to the
(irreducible part of) $\chi(\pi_1(\Sigma))$, on which $\hat\tau$ acts by the
identity,  $\hat\tau^* = id$.   When  $\alpha_{\Sigma}$ is $\su$-reducible, then
the cohomology groups are just the ordinary de Rham cohomology groups, tensored
with $su_2$, and the result is trivial.
\end{proof}

So in the case that  $\alpha_{\Sigma}$ is irreducible, the invariance of $L_t(B)$
under $\hat\tau^*$ follows directly from the claim.  In the case that 
$\alpha_\Sigma$ is $\su$-reducible, i.e.~when $t=0$, the Lagrangian subspace
$L_0(B)$ splits as the sum of
$r_B^*(H^1(B;\ada))$ and the anti-diagonal in $H^0(\Sigma;\ada) \oplus
H^0(\Sigma;\ada)$.  According to the claim, it is again the case that
$\hat\tau^*(L_0(B)) = L_0(B)$, and the theorem follows.
\end{proof}

The preceding argument contains most of the ingredients for comparing the $\Z$-grading, as
defined in~\cite{fs:graded}, for the groups $\HF^\mu(M)$ and  $\HF^\mu(\mtau)$.  The
spectral flow defined above is defined for paths of actual connections, rather than
gauge-equivalence classes of connections.  An fundamental observation is that it descends
to a function on
$\mathcal{A}/\mathcal{G}_0$, where $\mathcal{G}_0$ is the degree-$0$ gauge group. 
Likewise, the (perturbed) Chern-Simons function on connections descends to a function
$\widetilde{\CS}_h:\mathcal{A}/\mathcal{G}_0 \to \R$.  (Given a choice of trivial connection
$\Theta_M$, a lifting
$\tilde\alpha \in \mathcal{A}/\mathcal{G}_0$ of an \hflat connection, and a path of
connections from $\Theta_M$ to $\tilde\alpha$, the usual Chern-Weil integral over $M \times
\mathbf{I}$ defines $\widetilde{\CS}_h(\tilde\alpha)$.) 
Both $\SF$ and $\widetilde{\CS}_h$ depend on this choice of trivial connection.  

Fix any choice $\Theta_M$, which will
be used to pick out connections on all the other manifolds involved in the argument; for
starters the trivial connection on $\Sigma$ will simply be the restriction of $\Theta_M$.  
Recall construction~\ref{path} from the discussion leading up to the proof of
Theorem~\ref{sfequal}, in which we chose a specific path of connections $\alpha_t$ on $M$,
which were flat and irreducible on
$\Sigma$, with endpoint the trivial connection $\Theta_M$.  Simultaneously, we chose a
path of gauge transformations $g_t$ on $\Sigma$ with $g_t^*(\alpha_t) =
\tau^*(\alpha_t)$.   We make the convention in the subsequent discussion that any path of
connections on $M$ with endpoint
$\Theta_M$ should agree with this fixed model path near its endpoint.  Note that $g_0$ will
{\sl not} be the identity gauge transformation, because $g_t$ is of order $4$ for $t > 0$.

Choose a real number $\mu$ which is  not in the discrete set $ \{\widetilde{\CS}_h(\alpha)|\
\alpha
\in
\chi_h(M)\}$, and for each \hflat connection $\alpha \in \chi_h(M)$, pick a representative
$\tilde\alpha$ of its $\mathcal{G}_0$ orbit with
$\widetilde{\CS}_h(\tilde\alpha) \in (\mu,\mu+1)$.  The connections $\tilde\alpha$
form the basis of $\CF_*^\mu$, and the grading is defined in terms of spectral flow, where
one uses the same trivial connection, and path, as were used to define the Chern-Simons
invariant.  (We remind the reader that the spectral flow
$\SF(D_\alpha,D_\beta)$ changes by
$8\deg(g)$ when one replaces $\alpha$ by $g^*(\alpha)$.  Hence, if we consider specific
$\mathcal{G}_0$-representatives, the grading is actually $\Z$-valued, and not just
$\Z/8$-valued.)

\begin{theorem}\label{sfmu}
For any \hflat connection $\alpha$, the $\Z$-grading of
$\alpha^\tau$ is equal to the $\Z$-grading of $\alpha$.  Specifically, this means:
\begin{enumerate}
\item The \hflat  connection $\alpha^\tau$ constructed by
cutting and pasting satisfies $\widetilde{\CS}_h(\alpha) = \widetilde{\CS}_h(\alpha^\tau)$.
\item 
If $\mu \not\in \{\widetilde{\CS}_h(\alpha)|\ \alpha \in \chi_h(M)\}$, then   $\mu \not\in
\{\widetilde{\CS}_h(\beta)|\ \beta \in \chi_h(\mtau)\}$.   
\item $\SF_M(D_\Theta,D_\alpha) = \SF_{M^\tau}(D_\Theta,D_{\alpha^\tau})$
\end{enumerate}
\end{theorem}
\begin{proof}
Choose a path $\alpha_t$ of connections on $M$ from $\Theta_M$ to $\alpha$, agreeing with
the one from construction~\ref{path} near $t=0$.  As was remarked earlier, the whole path
can be extended to a path of connections $A_t$ on the basic cobordism $\wtau$.  The
restriction of $A_t$ to the boundary component $\mtau$ starts at a trivial connection
$\Theta_{\mtau}$ and ends at $\alpha^\tau$.  Because  $\alpha_{\Sigma,t}$ is flat, the
restriction (say, $\alpha_{S,t}$ of $A_t$ to the other boundary component $\mapt$ is flat. 
However, $\alpha_{S,0}$ is not the trivial connection, because the holonomy around $\zz$
is non-trivial. 

In a standard way, the family of connections $A_t$ fits together to give a connection on
the $5$-manifold $\wtau \times \mathbf{I}$; denote the restriction of this connection to
$\mapt \times \mathbf{I}$ by $A_S$.  Applying Stokes' theorem to the perturbed Chern-Weil
integrand 
$$\frac{1}{8\pi^2} \mathrm{Tr}((F_{A} + x) \wedge (F_A + x))$$ 
gives
\begin{multline*}\label{cstilde}
\widetilde{\CS}_M(\alpha) - \widetilde{\CS}_{\mtau}(\alpha^\tau) =
\frac{1}{8\pi^2}\int_{\mapt \times \mathbf{I}} \mathrm{Tr} (F_{A_S} \wedge F_{A_S})\\
+ \frac{1}{8\pi^2}\int_{\wtau} \mathrm{Tr} ((F_{A_1} + x) \wedge ( F_{A_1} + x))
- \frac{1}{8\pi^2}\int_{\wtau} \mathrm{Tr} ((F_{A_0} + x) \wedge (F_{A_0} + x))
\end{multline*}

Because each connection $\alpha_{S,t}$ is flat, the first integral on the right side
of this equation vanishes.  Similarly, since $A_0$ and $A_1$ are \hflat connections, the
other two integrals vanish, and the first item is proved.  The second one follows
immediately, and the third one is the content of Theorem~\ref{sfequal}.
\end{proof}

\subsection{Invariance of $HF_*$ under mutation}
We can now state and prove a more precise version of Theorem~\ref{hfmut}. 
\begin{theorem}\label{ytau} Let $M$ be an oriented homology $3$-sphere, with 
(instanton) Floer homology $HF_*(M)$, which contains a genus-$2$ surface,
and let $\mtau$ be the result of mutation along $\Sigma$.  Let $\ytau$ be the
orbifold constructed in section~\ref{extendsec}.  Then 
$$\ytau_*: HF_i(M) \xrightarrow{\cong}
HF_i(\mtau)
$$  
Similarly, if $\mu$ is a real number which is not the Chern-Simons
invariant of any flat connection on $M$, then $\ytau_*$ induces an isomorphism
on the $\Z$-graded instanton homology of Fintushel-Stern:
$$\ytau_*: HF^\mu_*(M)
 \xrightarrow{\cong} HF^\mu_*(\mtau).$$
\end{theorem}
\begin{proof} 
Choose a perturbation, as in Theorem~\ref{perturb}, so that the \hflat connections on
$M$ and $\mtau$ are all isolated smooth points.  This perturbation extends to a
`time-independent' perturbation of the ASD equations on $\ytau$, and so induces a map on
the $\Z$-graded theory as well as the usual Floer theory.  Choose a metric on $\ytau$ which
is generic, so that all moduli spaces on $\ytau$ are smooth of the correct dimension.  Let
$\alpha$ and $\beta$ be \hflat connections on $M$ and $\mtau$, respectively.   As
usual in Floer theory, the map induced by $\ytau$ is given on the chain level by
a matrix, whose $(\alpha,\beta)$ entry counts (with signs) the number
of points in the $0$-dimensional moduli space
$$ \M_h^w(\ytau;\alpha,\beta) $$
Since $\ytau$ is an orbifold, one needs to check that $\ytau_*$, as defined this
way, is in fact a chain map.  

As usual in Floer theory, this amounts to a constraint on the possible
non-compactness of $1$-dimensional moduli spaces
$\M_h^w(\ytau;\alpha,\beta)$.  One needs to see that there cannot be a
$1$-dimensional moduli space bubbling off at the singular points of the
orbifold.  Because $w_2$ is preserved in Uhlenbeck limits, the background
connection in such a circumstance would necessarily have non-trivial holonomy
at each cone point.  But (cf.~\cite{austin:so3,furuta:equivariant,furuta-hashimoto}) the minimum
dimension for a moduli space of ASD connections on $\Sigma(\RP^3)$ is $2$.
{}From the description of the cohomology of $\wtau$ in Lemma~\ref{whomology} ,
the rational cohomology of $\ytau$ is readily obtained.  Regarding $\ytau$ as
a rational homology manifold, we have
$$
b_2^+(\ytau) = 2, \ \chi(\ytau) = 4,\ \mathrm{and}\ \sigma(\ytau) = 0
$$
Since $b_2^+(\ytau) >1$, reducibles may be avoided, and so $\ytau_*$ yields a
well-defined chain map.  It remains to compute the degree of $\ytau_*$, and to
show that it induces an isomorphism on $\HF_*$ and $\HF_*^\mu$.

The degree is computed (compare~\cite{fs:pseudo}) by the
formula:
$$
\deg(\ytau_*) = 3(b_1(\ytau) - b_2^+(\ytau)) + \sum_{\mathrm{cone\ points}}
(d_\sigma + d_\chi)
$$
where the $d_\sigma$ and $d_\chi$ in the last sum refer to the signature and
Euler characteristic defects associated to each cone point.   Since the
holonomy of the bundle is non-trivial at each cone point, $d_\chi = 1$ and
$d_\sigma = 0$.  From the cohomology calculations, it follows that $
\deg(\ytau_*)  =0$.

The \hflat connections yield bases for the Floer chains on $M$ and $\mtau$,
and $\ytau_*$ is described as a matrix with respect to those bases. 
We will
show that $\ytau_*$, as a map on $\CF_*$, is represented by an upper
triangular matrix, with $\pm 1$ along the diagonal.  (The signs are presumably all $1$'s,
but we have not checked them.)
As explained in~\cite{fs:graded} (with an obvious adaptation to the orbifold setting), the
$\Z$-grading induces an increasing filtration on
$\CF_*$, and the chain map $\ytau_*$ respects the filtration; in matrix terms this means
that $\ytau_*$ is upper triangular.  Moreover, the proof of this statement (in Theorem 5.2
of ~\cite{fs:graded}) shows that if
$\alpha \in
\chi_h(M)$ and
$\beta\in \chi_h(\mtau)$ have the same grading, then any $A$ in a $0$-dimensional moduli
space $\M^w_h(\ytau)$ must in fact be $h$-flat.  (A similar argument,
in terms of an ordering on the \hflat connections, is presented
in~\cite[\S3.1]{braam-donaldson:knots}.)   In particular, the
$(\alpha,\alpha^\tau)$ entry is the number of \hflat extensions of $\alpha$ over $\ytau$,
which is one.  Lemma~\ref{smooth} says that this one point in
$\M^w_h(\ytau;\alpha,\alpha^\tau)$ is a smooth point, so it counts for $\pm 1$.   Thus
$\ytau_*$, which is a chain map, is an isomorphism, and therefore induces an isomorphism on
$\HF^\mu_*$ and $\HF_*$.

\end{proof}

\section{Mutation of $4$-manifolds}\label{4mfd}
The fact that the involution $\tau$ is in the center of the mapping class group
of the genus-$2$ surface leads to two types of cutting/pasting operations along certain
$3$-manifolds embedded in a $4$-dimensional manifold $X$.   Under the hypothesis that the
character varieties of the $3$-manifold in question is smooth, we will show that the
operation preserves the
$\su$-Donaldson invariant of
$X$. Both constructions have the additional feature that it
appears that gauge-theoretic invariants associated
to $\SU(3)$-bundles {\sl should} change.   There is no known construction
at this time for invariants associated to higher-rank bundles, unfortunately.

\subsection{Mutation along $3$-manifolds of Heegaard genus-$2$}
The following simple lemma is well-known: see~\cite{kobayashi:heegaard} and the
discussion
of problem 3.15 in~\cite{kirby:problems96}.
\begin{lemma}\label{heegaard}
Let $M^3$ be a closed $3$-manifold which admits a Heegaard splitting of 
genus $2$.  Then there is an involution $T:M \to M$ which preserves the
handlebodies and restricts to $\tau$ on the Heegaard surface $\Sigma$.
\end{lemma}
The statement of the lemma contains the proof-since $\tau$ commutes
with the attaching map, its 
extension over the genus-$2$ handlebodies fits together to define an involution
on $M$.

\begin{defn} Suppose that $X$ is a $4$-manifold, containing a  genus-two $3$-manifold
$M$.   Then
$X^\tau$, the mutation of
$X$ along
$M$ is defined to be the result of cutting and pasting $X$ along
$M$, using the involution $T$ described in the previous lemma.
\end{defn}
 
The free group on two letters has the same symmetry
property (Lemma~\ref{symmrep}) with regard to \su\ (or $\SO(3)$) representations
as does the fundamental group of a genus-$2$ surface.   This may be proven by
computing characters, as in~\cite{ruberman:mutation}, or by the following direct
 argument.  Write $A \in \su$ as $\exp(a)$, for $a \in {su}_2$. Now any element
$c \in  a^\bot \subset {su}_2$ of length $1/\pi$   will
have the property that 
$$ \exp(c)^{-1} A \exp(c) = A^{-1} = T_*(A).
$$
If $B$ is another element of \su, written as $\exp(b)$, then there is at least one
element $c$ in $a^\bot \cap b^\bot$ of length $1/\pi$.  Conjugating by $\exp(c)$ 
will take $A \to T_*(A)$ and $B \to T_*(B)$.  The same proof works for $\SO(3)$
representations.

\begin{theorem}\label{donaldson} Suppose that $M$ is a genus-$2$ homology sphere, for which
the \su\ character variety $\chi(M)$ is smooth. Then for any
Donaldson invariant $D_X$  which is defined on $X$, we have 
$$ D_X = D_{X^\tau}$$
\end{theorem}
\begin{remark}
As usual in gauge theory, the hypothesis of smoothness means that $\chi(M)$ is a smooth
manifold, whose tangent space at $\alpha$ is given by $H^1(M;\ada)$.
\end{remark}
\begin{proof}
Under the smoothness hypothesis, it is known
(compare~\cite{gompf-mrowka,mmr}) that one can express the Donaldson invariant in terms of
ASD connections on the components of $X - M$, exponentially decaying to flat connections on
$M$.  But the symmetry property described above shows that $T_*$ acts trivially on the 
character variety $\chi(M)$, and hence does not change the gluing picture, so the
Donaldson invariant does not change.   
\end{proof}

It seems reasonable that the theorem should continue to hold in the general case, when
$\chi(M)$ is not discrete.  One would need to first find a $T$-equivariant perturbation $h$
of the Chern-Simons invariant, and then to show that the induced action of $T$ on the
$h$-flat connections is the identity.   It is not hard to find an equivariant perturbation,
by arranging that the link $L$ along which the perturbation is supported to
be $T$-invariant.  Unfortunately, we have been unable to prove that this means that the
action of $T^*$ on $\chi_h(M)$ is the identity.  

If there were a Donaldson-type theory for the group $\SU(3)$, then mutant
$4$-manifolds would not likely have the same invariants.  Here is an example, based on
calculations of Hans Boden, which is waiting for an appropriate theory.
\begin{example}
Suppose that $M = \Sigma(a_1,a_2,a_3)$ is a Seifert-fibered space with $3$
exceptional fibers.  Then $M$ has a genus-$2$ Heegaard splitting.  View $M$  as
$S^1
\times (S^2 - 3\ \rm{int}D^2) \cup_3 S^1 \times   D^2$, where the gluing maps give
rise to the multiplicity of the fibers.  One of the genus-$2$ handlebodies is
given by the union of two of the $S^1 \times D^2$, together with a thickened arc
joining them.  (That the complement is also a handlebody is left to the reader to
verify!).  In~\cite[\S5.1]{boden:brieskorn},   Boden shows that for $\Sigma(2,3,7)$ there
are $4$ irreducible $\SU(3)$ representations, two of which are the complexification of an
$\SO(3)$ representation, and which are therefore invariant under $T_*$.  The other two,
denoted $\rho_3$ and $\rho_4$, are interchanged by $T_*$.  One would expect that $T_*$
would thus act non-trivially on an
$\SU(3)$ Floer-type theory.  A similar phenomenon occurs for $\Sigma(2,3,6k\pm 1)$ and
presumably most Brieskorn spheres.

%%The standard presentation for the fundamental group of $M$ is
%%$$\langle x_1,x_2,x_3, h | h\ \mathrm{central},\ x_i^{a_i} = h^{-b_i},\ x_1x_2x_3
%%= h^{-b_0} \rangle
%%$$
%for appropriate values of the $b_i$.
\end{example}

\begin{remark}
If the Brieskorn sphere
$\Sigma$ is viewed as the link of a surface singularity (defined over the reals), then $T$
may be identified with the involution coming from complex conjugation.  In particular, $T$
extends over such $4$-manifolds as the Milnor fiber and canonical resolution of the
singularity.  So, for example, mutation of a $K3$ surface along $\Sigma(2,3,7)$ does not
produce any new $4$-manifolds.  As pointed out by Tom Mrowka, $T_*$ does not act
trivially on the Floer-theory associated to the Seiberg-Witten equations.  For it may
be seen from the description of $T$ as complex conjugation and the computations
in~\cite{moy} that $T_*$ in fact acts as the well-known involution in the Seiberg-Witten
theory.
\end{remark}

\subsection{Mutation along genus-$2$ mapping tori}
Suppose that  $\varphi$ is a diffeomorphism of $\Sigma$, and form the mapping torus
$\mapphi$.  The fact that $\tau$ is in the center of the genus-$2$ mapping class group
implies that there is an involution $T:\mapphi \to \mapphi$ whose restriction to a  fiber
is $\tau$.   If $\mapphi$ is embedded in a $4$-manifold $X$, then as above, we define the
mutation of $X$ along $\mapphi$ by cutting and pasting via $T$.  In order to carry out the
argument described above, it is important to avoid the reducible flat
connections on $\mapphi$.  One way to do this is to assume that the Donaldson invariant
being computed is associated to an $\SO(3)$ bundle with $\langle w_2, \Sigma \rangle \neq
0$.   We need the analogue of Lemma~\ref{symmrep}:
\begin{lemma}\label{maptrep}
Let $T$ be the involution on $\mapphi$ induced by $\tau$, and let $\rho$ be an
$\SO(3)$ representation of $\pi_1(\mapphi)$ such that $\langle w_2(\rho), \Sigma \rangle
\neq 0$.  Then $\rho\circ T_*$ is conjugate
to $\rho$.
\end{lemma}
\begin{proof}
The fundamental group of $\mapphi$ has a standard presentation as
$$
\langle t, \pi_1(\Sigma) | t^{-1}gt = \varphi_*(g)\ \forall g\in \pi_1(\Sigma) \rangle
$$
The main point to notice is that if $\rho$ is any representation (\su \ or $\SO(3)$), then
$\rho(t)$ is determined by $\rho_\Sigma$, up to elements in the centralizer of
$\rho(\pi_1(\Sigma))$.  Let $\gamma$ be any element with $\gamma^{-1} \rho_\Sigma \gamma =
\rho_\Sigma$.  Then it is easy to check that $\gamma^{-1} \rho(t) \gamma = h \rho(t)$ for
some $h$ in the centralizer of $\rho(\pi_1(\Sigma))$.  In order to prove the lemma we must
show that in fact $h$ is trivial.  Under the hypothesis on
$w_2$, there are only two possibilities: either the centralizer is trivial (in which case
the lemma follows directly), or $\rho(\pi_1(\Sigma))$ is a $\Z_2 \oplus \Z_2$ subgroup,
which is its own centralizer.   But since every element of  $\Z_2 \oplus \Z_2$ has order
$2$, one can take $\gamma$ to be trivial, in which case it certainly conjugates $\rho(t)$
to itself.
\end{proof}

With this lemma in hand, the proof of Theorem~\ref{donaldson} yields an analogous theorem
about mutation along mapping tori.
\begin{theorem}
Suppose the mapping torus $\mapphi$ is embedded in $X$, and that $D^w$ is a Donaldson
invariant associated to a bundle with $\langle w_2, \Sigma \rangle \neq 0$.  Then $D^w(X) =
D^w(X^T)$.
\end{theorem}

\end{document}